\documentclass[11pt]{amsart}
\usepackage{amssymb,latexsym,amsmath}
\input{xypic}
\xyoption{all}
\numberwithin{equation}{section}

\newcommand{\R}{{\mathbb R}}

\numberwithin{equation}{section}
\newtheorem{theorem}[equation]{Theorem}

\begin{document}

\title[]{A simple generalised Plancherel formula for compactly induced characters}
\date{\today}
\author{\bf Chaitanya Ambi}
\address{Indian Institute of Science Education and Research, Dr.\,Homi Bhabha Road, Pashan, Pune 411008,  INDIA.}
\email{chaitanya.ambi@students.iiserpune.ac.in}

\maketitle

\begin{abstract} 
The aim of this article is to present a simple generalized Plancherel formula for a locally compact unimodular topological group $G$ of type I. This formula applies to the functions representing c-$ Ind_{U}^{G}\psi $ for a unitary character $ \psi $ of a closed unimodular subgroup $ U $ of $ G $. This specializes to the Whittaker-Plancherel formula for a split reductive $p$-adic group of Sakellaridis-Venkatesh and differs from that of a quasi-split $p$-adic group due to Delorme. Furthermore, it also applies to certain metaplectic groups  and other interesting situations where the local theory of distinguished representations has been studied. 
\end{abstract}

\tableofcontents

\section{\bf Introduction and statement of the main theorem}
Let $ G $ be a unimodular, locally compact topological group of type I with the unitary dual $ \hat{G} $.  Fix a Haar measure on $ G $. Let $ \mathcal{L}^{2}(G) $ denote the space of complex-valued functions on $ G $ which are square integrable with respect to the left Haar measure. The left regular representaion of $ G $ on $ \mathcal{L}^{2}(G) $ can be disintegrated into irreducible representations. This permits one to express a function within a suitable class in $ \mathcal{L}^{2}(G) $ in terms of functions on $ \hat{G} $ via the Fourier transform. The Plancherel formula  reflects the structure of $ \hat{G} $. \\

For a closed subgroup $ U $ of $ G $ having a continuous unitary character $ \psi $, a function in the space representing the compact induction c-$ Ind_{U}^{G}\psi $ also affords a generalised Plancherel formula which provides further information about $ \hat{G} $. Such formulae have been worked out in various settings, especially when $ G $ is a $ p $-adic group or a Lie group with certain properties and $ U $ is the unipotent radical of a Borel subgroup of $G $ (notably, the Whittaker-Plancherel formulae by Baruch and Mao [1], Delorme [3] and also Sakellaridis and Venkatesh [10]).\\
      In this article, we derive a generalised  Plancherel formula for a broad class of groups  which generalises  such Whittaker-Plancherel type formulae in the pointwise case. We prove the existence of a locally integrable kernel for the  distribution character and  compute it in terms of the kernel of the ordinary Plancherel formula for $ G $ whenever the latter exists. We also establish the absolute continuity of the corresponding Plancherel measure
      with respect to that on $ \hat{G} $. Thus, our result presented below facilitates an explicit derivation of a generalised Plancherel formula in several situations including  those where the  formulae known so far do not apply (such as certain metaplectic groups).\\

      Assume that $ U $ is unimodular. For $ f\in \mathcal{C}_{c}(G) $ and  $g\in \mathcal{L}^{\infty}(U)$, set \begin{equation}(g*_{U}f)(x)=\int_{U}^{}g(u)f(u^{-1}x)du.\end{equation}  
      
  Let $ \mathcal{C}_{c}(G) $ denote the space of  compactly supported, complex-valued continuous functions on $ G $. For $ f \in  \mathcal{C}_{c}(G) $, define the operator 
\begin{equation}
\hat{f}(\pi) = \int_{G}^{}f(x)\pi(x)dx,
\end{equation}  which is known to be compact and Hilbert-Schmidt (hence of the trace class).\\
 Since we have fixed a Haar measure on $ G $, there exists a unique Plancherel  measure $\mu_{\pi}$ on $ \hat{G} $ such that the following pointwise Inversion Formula holds (see [4]):
\begin{equation}
h(1)= \int_{\hat{G}^{}}\Theta_{\pi}(h)d\mu_{\pi}  \\, \quad \forall  h \in \mathcal{C}_{c}(G),
\end{equation} where the distribution $ \Theta_{\pi} $ is defined as 
\begin{equation}
\Theta_{\pi}(h)= Tr[\hat{h}(\pi)].
\end{equation}

For a continuous unitary character $ \psi $ of $ U $, define the space 
\begin{equation}
\textit{c-Ind }_{U}^{G}\psi := \{ W:G \rightarrow  \mathbb{C} :\text{ W is compactly supported modulo U 
	 and satisfies }  W(ug)=\psi(u) W(g)\},
\end{equation} ($ u \in U $).
This space serves as a model for compact induction of $ \psi $. Define 
\begin{equation}
(\mathcal{W}_{\psi}f)(g) = (\psi*_{U}f)(g). 
\end{equation}
It can be shown (see [4]) that the map $ f \rightarrow \mathcal{W}_{\psi}f $ is surjective.\\
In many situations, there exists (except for a set of measure zero in $ \hat{G} $) a locally integrable function $ \theta_{\pi}:G\rightarrow \mathbb{C}$  $(\pi \in  \hat{G}) $ which is  constant on conjugacy classes of $ G $ and satisfies
\begin{equation}
\Theta_{\pi}(f)=\int_{G}^{}f(g)\theta_{\pi}(g)dg, \space (f \in \mathcal{C}_{c}(G)).
\end{equation}  This occurs in the case of real reductive groups as well as  reductive $ p $-adic groups (see  [6],[7] and [8]). The existence of $\theta_{\pi}  $, which we shall term as the \textit{ kernel} of $ \Theta_{\pi} $, permits us to derive a simple generalised Plancherel formula as follows:

\begin{theorem}[Generalised Plancherel formula for a type-I group]
Assume that the distribution character $ \Theta_{\pi} $ of $ G $ has a locally integrable function $ \theta_{\pi} $ as its kernel. Then, there exists a distribution $ \Phi_{\pi}^{\psi} $ (defined except on a set of measure zero for $ \pi \in \hat{G} $) such that the function $ \mathcal{W}_{\psi}f \in \textit{c-Ind}_{U}^{G}\psi $  satisfies
\begin{equation}
(\mathcal{W}_{\psi}f)(1)= \int_{\hat{G}}^{}\Phi_{\pi}^{\psi}(f)d\mu_{\pi}.
\end{equation}  
The distribution $ \Phi_{\pi}^{\psi} $ is given explicitly by
\begin{equation}
\Phi_{\pi}^{\psi}(f) = \int_{G}^{}f(g)(\bar{\psi}*_{U}\theta_{\pi})(g)dg,
\end{equation}
 where kernel $  \bar{\psi}*_{U}\theta_{\pi} $ of $ \Phi_{\pi}^{\psi} $ exists as a locally integrable function on $ G $.

\end{theorem}
(See also the Remark in Section (3.1)).\\

  \textbf{Acknowledements:} The author is supported by a grant by NBHM. The author would also like to thank Prof. A. Raghuram for his suggestions, support and encouragement.
  
\bigskip
\section{\bf Proof of the generalised Plancherel theorem for a type-I unimodular group}
We shall maintain the notation as in the Introduction.\\
   Since $ U $ is closed in $ G $, it is also locally compact and second countable. Hence, $U$ is $ \sigma $-compact and we can select a sequence $ K_{n} $ of symmetric compact subsets of $ U $ such that 
\begin{equation}
K_{i}\subset K_{i+1} \quad \forall i \in \mathbb{N} \text{ and } \cup K_{i}=U.
\end{equation}
Define a sequence of  bounded, compactly supported functions $ {\psi_{n}: n \in \mathbb{N}} $ by 
\begin{equation}
\psi_{n}(u)= \psi(u) \mathbf{1}_{K_{n}}(u) \quad(u \in U)
\end{equation}
For each function  $ f \in \mathcal{C}_{c}(G) $ and $ x,y \in G $, we have the estimate 
\begin{equation}
|\psi_{n}*_{U}f(x) - \psi_{n}*_{U}f(y)| \le ||\psi_{n}||_{\mathcal{L}^{\infty}(U)} \int_{U}^{}|f(u^{-1}x)-f(u^{-1}y)|du.
\end{equation}
Noting the uniform continuity of $f$ on its support and that no $ \psi_{n} $ exceeds 1 in magnitude, we conclude that each $ \psi_{n}*_{U}f $ is continuous on $ G $. Further, both  $ \psi_{n} $ and $f $ are compactly supported and hence $ \psi_{n}*_{U}f \in \mathcal{C}_{c}(G) $ for each $ n $. Hence, the formula (2.2) applies to $ \psi_{n}*_{U}f $ and we have 
\begin{equation}
(\psi_{n}*_{U}f)(1)= \int_{\hat{G}}^{}\Theta_{\pi}(\psi_{n}*_{U}f)d\mu_{\pi}.
\end{equation}
\begin{equation}
=\int_{\hat{G}}^{} \int_{G}^{}\theta_{\pi}(g)(\psi_{n}*_{U}f)(g)dgd\mu_{\pi}.
\end{equation}
\begin{equation}
=\int_{\hat{G}}^{} \int_{G}^{}\int_{U}^{}\theta_{\pi}(g)f(u^{-1}g)\psi_{n}(u)dudgd\mu_{\pi}.
\end{equation}
Now, we use Fubini's Theorem to interchange the inner two integrals. This stands justified as the whole integrand is integrable with respect to the product measure $ dgdu $ (as $ f $ is compactly supported) and both $ dg $ and $ du $ are $ \sigma $-finite (owing to the local compactness and second countability of $ G $ and $ U $, respectively). 
\begin{equation}
=\int_{\hat{G}}^{} \int_{U}^{}\psi_{n}(u)(\int_{G}^{}\theta_{\pi}(g)f(u^{-1}g)dg)dud\mu_{\pi}.
\end{equation}
Substituting $ g=ux $ in the innermost Haar integral, we get
\begin{equation}
=\int_{\hat{G}}^{} \int_{U}^{}(\int_{G}^{}\theta_{\pi}(ux)f(x)dx)\psi_{n}(u)dud\mu_{\pi}.
\end{equation}
Since we have not changed the integrand, Fubini's Theorem applies as before.
\begin{equation}
=\int_{\hat{G}}^{} \int_{G}^{}f(x)(\int_{U}^{}\theta_{\pi}(ux)\psi_{n}(u)du)dxd\mu_{\pi}.
\end{equation} 
Noting the unimodularity of $ U $ and that $ \bar{\psi_{n}}(u)=\psi_{n}(u^{-1}) $,  (bar denotes complex conjugation throughout this article), we obtain
\begin{equation}
=\int_{\hat{G}}^{} \int_{G}^{}f(x)(\int_{U}^{}\bar{\psi_{n}}(u)\theta_{\pi}(u^{-1}x)du)dxd\mu_{\pi},
\end{equation}
We recognise the above expression as 
\begin{equation}
=\int_{\hat{G}}^{} \int_{G}^{}f(x)(\bar{\psi_{n}}*_{U}\theta_{\pi})(x)dxd\mu_{\pi}.
\end{equation}
Note that $ \bar{\psi_{n}}*_{U}\theta_{\pi} $ exists as a locally integrable function on $ G $ because $ \psi_{n}$ is compactly supported and $\theta_{\pi} $ is locally integrable (see Prop. 5.4.25 on p.283 in [5]).\\
Now, we restrict ourselves to an arbitrary compact subset $ K $ of positive measure in $ G $. By Prop. 4.4.24 on p.282 in [5], $ \bar{\psi}*_{U}\theta_{\pi}  $ exists as an integrable function on $ K $.\\
It remains to  establish the convergence of $ \bar{\psi_{n}}*_{U}\theta_{\pi} $ to $ \bar{\psi}*_{U}\theta_{\pi}  $ in $ \mathcal{L}^{1}(K) $. This follows from the Dominated Convergence Theorem once  we observe that $ \bar{\psi_{n}} $ converges to $ \psi $ pointwise and that the integrable function $ 1+|\theta_{\pi}| $ dominates each $ \bar{\psi_{n}}*_{U}\theta_{\pi} $ almost everywhere on $ K $. This proves the existence of $ \bar{\psi}*_{U}\theta_{\pi} $ as a locally integrable function on $ G $, which  we have denoted  by $ \phi_{\pi}^{\psi} $ in the statement of Theorem (1.8). 
If we denote the distribution with $ \phi_{\pi}^{\psi} $ as its kernel by $ \Phi_{\pi}^{\psi} $, we have 
\begin{equation}
(\mathcal{W_{\psi}}f)(1)=\int_{\hat{G}}^{}\Phi_{\pi}^{\psi}(f)d\mu_{\pi}\quad ( f \in \mathcal{C}_{c}(G)),
\end{equation}
which is the pointwise generalised Plancherel formula for $ \mathcal{W_{\psi}}f $. If we define the convolution of $ \psi $ with a distribution in a similar manner, $ \Phi_{\pi}^{\psi} $ equals
\begin{equation}
\Phi_{\pi}^{\psi} = \psi*_{U}\Theta_{\pi}.
\end{equation}
This completes the proof of Theorem (1.8) once we note that the Plancherel measure $ d\mu_{\pi} $ remains unaltered throughout the computation.
\bigskip
\section{\bf Comparison with other known Whittaker-Plancherel formulae} 
It is well-known that  a real or $ p $-adic reductive group  $G$ is of type I (see [2]). The unipotent radical $U$ of its standard Borel subgroup is unimodular as well as closed. Hence, $U$ is also locally compact as well as second countable. Theorem (1.8) is applicable and  provides a broader framework for the following Whittaker-Plancherel formulae : 

\subsection{\bf Sakellaridis-Venkatesh}
Sakellaridis and Venkatesh [10] have obtained a Whittaker-Plancherel formula (more precisely, the Parseval's Identity of isometry) for Whittaker functions on a split reductive algebraic group. They also prove the absolute continuity of the Whittaker-Plancherel measure with respect to the Plancherel measure for the group. Theorem (1.8) not only improves the formula of Sakellaridis and Venkatesh to a Whittaker-Plancherel  pointwise inversion formula, but also allows us to find the Radon-Nikodym derivative of the involved measure, provided the measure is suitably normalised (namely, equal to 1).\\

\textbf{Remark:}
Since  Plancherel measures are defined only upto absolute continuity, we could have adopted the alternative viewpoint that the Plancherel measure changes to $ \nu_{\pi}^{\psi}$ while the distribution $ \Theta_{\pi}  $ remains unchanged in the Whittaker-Plancherel case, i.e.,
\begin{equation}
\Theta_{\pi}d\nu_{\pi}^{\psi}=\Phi_{\pi}^{\psi}d\mu_{\pi}.
\end{equation}
We conjecture that the Radon-Nikodym derivative would then be equal to 
\begin{equation}
 d\nu_{\pi}^{\psi}/d\mu_{\pi} = mult(\pi, \textit{c-Ind}_{U}^{G}\psi),
\end{equation}\textit{where $ mult(\pi, Ind_{U}^{G}\psi) $ denotes the multiplicity of $ \pi $ in
 c-$ Ind_{U}^{G}\psi $.}\\
The multiplicity is known to be finite in several cases, such as that of degenerate Whittaker models (see [9]).

\subsection{\bf Delorme}
Delorme [3] has obained an inversion formula for a Whittaker function using the matrix coefficients of certain representations induced parabolically and a version of the Schur Orthogonality relation for such coefficients. Assuming a character $ \psi $ to be nondegenerate, Delorme's formula expresses a $ \psi $-Whittaker function in terms of certain transforms of the function itself. Theorem (1.8) expresses the value of a $ \psi $-Whittaker function at unity in terms of the compactly supported function which generates it and the kernel  $ \theta_{\pi} $ of a representation $\pi \in \hat{G}$. This is substantially different from the formula obtained by Delorme.

\subsection{\bf Baruch-Mao}
 Baruch and Mao [1] have obtained a Whittaker-Plancherel formula for $SL(2,\mathbb{R})$ by means of its characters expressed explicitly in terms of Bessel functions. Since $SL(2,\mathbb{R})$ is a connected semisimple Lie group, Theorem (1.8) evidently applies to it.
\bigskip

\subsection{\bf The Whittaker-Plancherel formula for metaplectic groups}
We now discuss the applicability of Theorem (1.8) to  metaplectic groups. Let $Sp_{2n}(\mathbb{F}) $ be the symplectic group over a local field $ \mathbb{F} $ of characteristic zero ($ n \in \mathbb{N} $). Consider a group $ Mp_{2n}(\mathbb{F}) $ which is a finite cover of $ Sp_{2n}(\mathbb{F}) $ (and thus is second countable and locally compact). Further, $ Mp_{2n}(\mathbb{F}) $ is semisimple and hence of type I. We conclude this article by mentioning that Theorem (1.8) is applicable to $ Mp_{2n}(\mathbb{F}) $ as well.

\end{document}